\documentclass[]{OAGM}
\OAGMarXiv{1404.3538}

\usepackage{graphicx}

\title{Code Minimization for Fringe Projection Based 3D Stereo Sensors by Calibration Improvement}

\author{Christian Br\"auer-Burchardt, Peter K\"uhmstedt and Gunther Notni\\
  Fraunhofer Institute Applied Optics and Precision Engineering Jena, Germany}

\begin{document}
\maketitle

\begin{abstract}
Code minimization provides a speed-up of the processing time of fringe projection based stereo sensors and possibly makes them real-time applicable. This paper reports a methodology which enables such sensors to completely omit Gray code or other additional code. Only a sequence of sinusoidal images is necessary. The code reduction is achieved by involvement of the projection unit into the measurement, double triangulation, and a precise projector calibration or significant projector calibration improvement, respectively.
\end{abstract}

\section{Introduction}

Precise surface measurements become more and more a matter of fringe projection profilometry. However, the requirements concerning spatial resolution, accuracy, and especially measurement speed are always increasing. Whereas faster and better hardware components fulfill these requirements on the one hand, better algorithms may improve measurement systems without considerable higher costs.

High-speed applications of precise contactless surface measurements which do not interrupt production process are of special interest in the field of industrial quality check and control. Measurement systems based on fringe projection technique are usually limited: projection and observation rate using common hardware components is typically not higher than $180\,Hz$, and long projection sequences decrease the 3D measurement frequency.

Typically, the projection sequence consists of a number of phase shifted sinusoidal images and an additional (e.g. Gray code \cite{Sansoni-1999-AO}) sequence in order to ensure uniqueness of the phase values, the so called phase unwrapping \cite{Zhang-1999-AO}. An overview over existing fringe projection methods is given by Zhang \cite{Zhang-2010-OLE}.

However, in order to realize real-time processing, additional code for phase unwrapping should be omitted. Methods for like Multi-frequency-technique \cite{Li-2005-OE} are out of interest because they require additional code. Two alternatives are advertised: either continuous phase tracking, as e.g. proposed by Herr\'{a}ez \cite{Herraez-2002-AO} or geometric phase unwrapping using additional geometric information, e.g. concerning the measurement volume (see \cite{CBB-2013-SPIE}).

Indeed, Herr\'{a}ez' algorithm has two main disadvantages. First, it has relative low robustness in case of objects with sharp edges and hidden parts leading to phase jumps. Second, calculation effort is quite extensive. Other interesting approaches for phase unwrapping algorithms using multiple projector or camera positions are presented by Ishiyama et al. \cite{Ishiyama-2007-AO}, or Young et al. \cite{Young-2007-CVPR}. 

In this paper, a new method will be described which realizes 3D point determination without additional code. It is achieved by phase unwrapping using double triangulation and evaluation of multiple 3D point candidates.

\section{Measurement Principle}

Using structured light for finding point correspondences leads to the principle of phasogrammetry which is the combination of photogrammetry and fringe projection in closed mathematical form (see \cite{Schreiber-2000-OE}).
A projection unit projects a well-defined sequence of fringe images onto the measurement object which is observed by one or more cameras $C_{i}$. One sequence of projected fringe images usually consists of a Gray code sequence (see e.g. \cite{Sansoni-1999-AO}) which realizes the uniqueness of the fringes and a sequence of sinusoidal fringe patterns.

The observed fringe image sequences are processed resulting in phase images for each measurement position and camera. After rotation of the fringe pattern by $90^{\circ}$, the sequence may be projected and observed again resulting in a second phase image for each camera. The phase values correspond to image coordinates in the projector image plane. The resulting 3D points are obtained by triangulation between the coordinates of the camera and the projector (CP-mode), or between corresponding points of two cameras (CC-mode). Triangulation can be regarded as standard procedure in photogrammetry (see \cite{Luhmann-2006-Wiley}).

\begin{figure}[ht]
  \centering\includegraphics[width=0.70\textwidth]{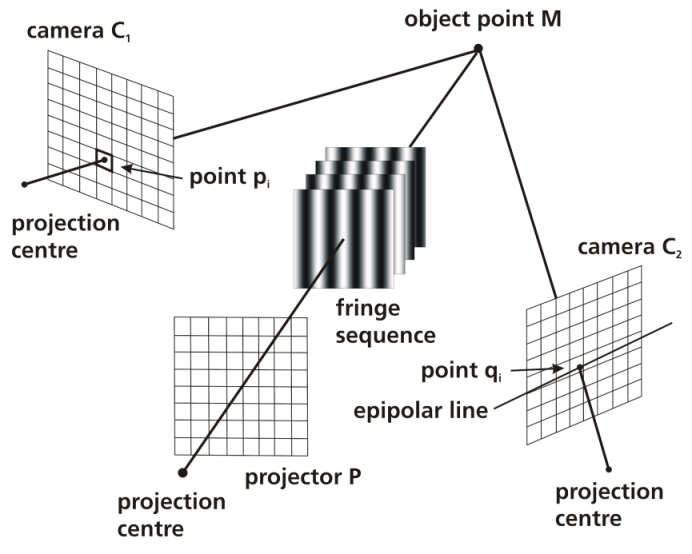}\\
  \caption{Principle of phasogrammetry using epipolar geometry.}
  \label{fig:fig1}
\end{figure}

Epipolar geometry is a well-known principle which is often used in photogrammetry when stereo systems are present as for example described in \cite{Luhmann-2006-Wiley}. It is characterized by an arrangement of two cameras observing almost the same object scene. A measurement object point $M$ defines together with the projection centres $O_{1}$ and $O_{2}$ of the cameras a plane $E$ in the 3D space (see also Fig.~\ref{fig:fig1}). The images of $E$ are corresponding epipolar lines concerning $M$. When the image point $p_{i}$ of $M$ is selected in camera $C_{1}$ image, the corresponding point $q_{i}$ in camera $C_{2}$ image must lie on the corresponding epipolar line. This restricts the search area in the task of finding corresponding points. 


\section{The New Approach}

\subsection{Finding Correspondences and 3D Point Calculation}

First the concept of the new method should be briefly outlined. It is based on the double application of triangulation procedure in order to calculate the resulting 3D points. The projected and observed sequence of sinusoidal fringe images is used to calculate one image of wrapped phase values for each camera. The principle of epipolar geometry should be applied in order to support finding phase correspondences.

\begin{figure}
		\centering
			\includegraphics[width=0.99\textwidth]{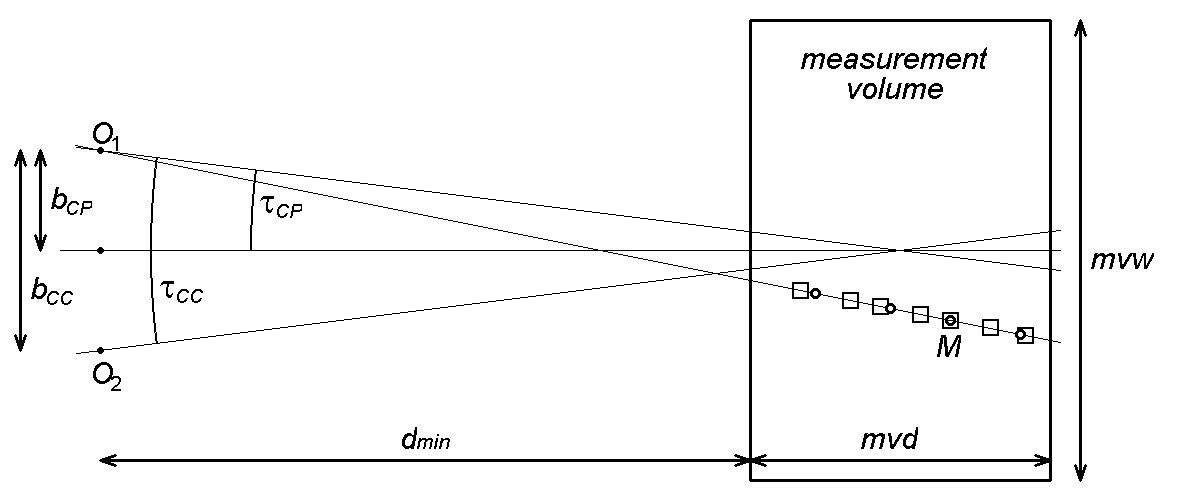}
					\caption{Geometric situation of the sensor arrangement (view from above), boxes represent CC-candidates and circles represent CP-candidates of the object point $M$ measurement.}
	\label{fig:fig2}
	\end{figure}

The idea is the following. For each camera pixel of camera $C_{1}$ all potential sub-pixel exact corresponding candidates in the image of camera $C_{2}$ are found and converted to potential 3D point candidates by triangulation (CC-mode). The obtained list is compared to those candidates found by application of the CP-mode using camera $C_{1}$ and the projector for triangulation (see Fig.~\ref{fig:fig2}).

Let a fringe projection based stereo sensor be given. The kind of the measurement objects, the arrangement of the optical components, and the adjustment of the lenses (range of sharp observation) restrict the measurement volume (MV) in its depth ($mvd$). This restriction leads to reduction of the search area for the corresponding phase value $\phi$ onto a segment $s$ in the image of camera $C_{2}$ using the CC-mode. All positions on $s$ corresponding to the phase value $\phi$ are candidates for point correspondence. The number $n_{CC}$ of such candidates depends on the measurement volume width ($mvw$) and depth ($mvd$), the number $N$ of projected fringes (or equivalently the fringe period length $\lambda_{N}$) and the main triangulation angle $\tau_{CC}$ (or equivalently the length $b_{CC}$ of the baseline between the projection centers $O_{2}$ and $O_{2}$  of the two cameras and the minimal object distance $d_{min}$). More details are given in \cite{CBB-2013-SPIE}. For illustration see Fig.~\ref{fig:fig2}. 

The relationships between the parameters can be expressed by the formula:

\begin{equation}
n = mvd \cdot \frac{N \cdot b - mvw}{mvw \cdot d_{min}}
\end{equation}

Certainly, $mvd$ corresponds to the length of the illuminated area by the $N$ projected fringes at the distance $d_{min}$. Because $n$ should be a natural number, the $n$ obtained by (1) must be rounded up.

	
Applying the CP-mode the same holds. However, the number $n_{CP}$ of candidates will be approximately only the half, if the projector is symmetrically arranged between the cameras (see Fig.~\ref{fig:fig2}). For each candidate the resulting 3D point candidate is generated by triangulation for both candidate lists (CC-mode and CP-mode) and all 3D points of CC are compared to those obtained by CP-mode. 

\subsection{Rejection of False Candidates}

All pairs of 3D points with a 3D distance $dst$ bigger than a given threshold $thr$ are rejected. This threshold $thr$ can be found by analysis of the typical phase error or can also be determined by experiments.
The ratio between false positive 3D points and correctly determined should be under a certain level which makes it easy to identify them as outliers in the resulting 3D point cloud and remove them manually or automatically.

\subsection{Outlier Removal}

There are two possible ways for outlier removal depending on the decision strategy for correct candidate selection. The first strategy (mode $m_{1}$) is to select exactly one candidate using the minimal 3D distance between the CC- and the CP-candidate. In this case the false candidate should be removed from the 3D point cloud. The second strategy (mode $m_{2}$) selects all candidates with distances $dst$ below $thr$ and only candidates with a distance bigger than $thr$ are removed from the 3D point cloud.

Which strategy will be selected should be decided after first corresponding experiments. In the current state of our algorithm the outlier removal was performed manually in the 3D point cloud (see next section). Alternatively, an automatic removal using point clustering and probabilistic properties of the produced clusters may be realized (see \cite{CBB-2011-Constraints}).

\subsection{Calibration Improvement by Precise Projector Distortion Calculation}

A careful calibration of fringe projection based 3D stereo sensors typically provides sufficient measurement accuracy according to the requirements. Using the CC-mode for triangulation, the calibration parameters of the projection unit are not necessary for the measurement process. Consequently, projection unit is often not or only insufficiently calibrated. However, an accurate projector calibration is necessary for our new code minimized method.

Especially the distortion is usually not sufficiently determined. This is often due to the fact, that projector distortion cannot be well described by a distortion function, because of the location of the principal point outside the projector chip. Here a very careful determination of the projector distortion must be performed. One possibility to realize this will be described by the new three step method presented in the following.

The precondition for the new method of projector distortion determination is an accurate calibration of the two cameras including camera distortion. There are several ways to evaluate the quality of the current calibration. One possibility is to check the measured surface of a well-known object concerning lengths and flatness of plane surface parts (as e.g. proposed in \cite{VDI-2008-V}), and to analyze epipolar line error and scaling error (see \cite{CBB-2011-ICIAP}).

A second precondition is an initial calibration of the projectiion unit.

In a first step we improved the intrinsic and extrinsic camera parameters of the two cameras of the sensor using the method proposed in \cite{CBB-2011-ICIAP} by minimization of the epipolar line error $\Delta E$ and the scaling error $\Delta S$. In the next step the intrinsic and extrinsic calibration parameters of the projector were improved using the same method.

In the third step the distortion matrix of the projector was improved the following way. A plane surface was placed into the measurement volume at three distances: at $d_{min}$, centrally, and at $d_{min} + mvd$. Measurements using CC-mode were performed. A set $P_{i}$ of 3D points was obtained for each measurement $i (i=1,2,3)$, assumed to be correct. Then a CP-triangulation was performed obtaining three sets $P_{i}'$ of corresponding 3D points. The elements of the $P_{i}'$ are assumed to be distorted. Now, every pair of corresponding 3D points determines a unique correction vector $d = (dx, dy)$ obtained by back-projection of the $P_{i}$ into the projector image plane $I_{p}$. The vectors $d_{i}$ describe the residuals in $I_{p}$.

\begin{figure}[ht]
  \centering\includegraphics[width=1.0\textwidth]{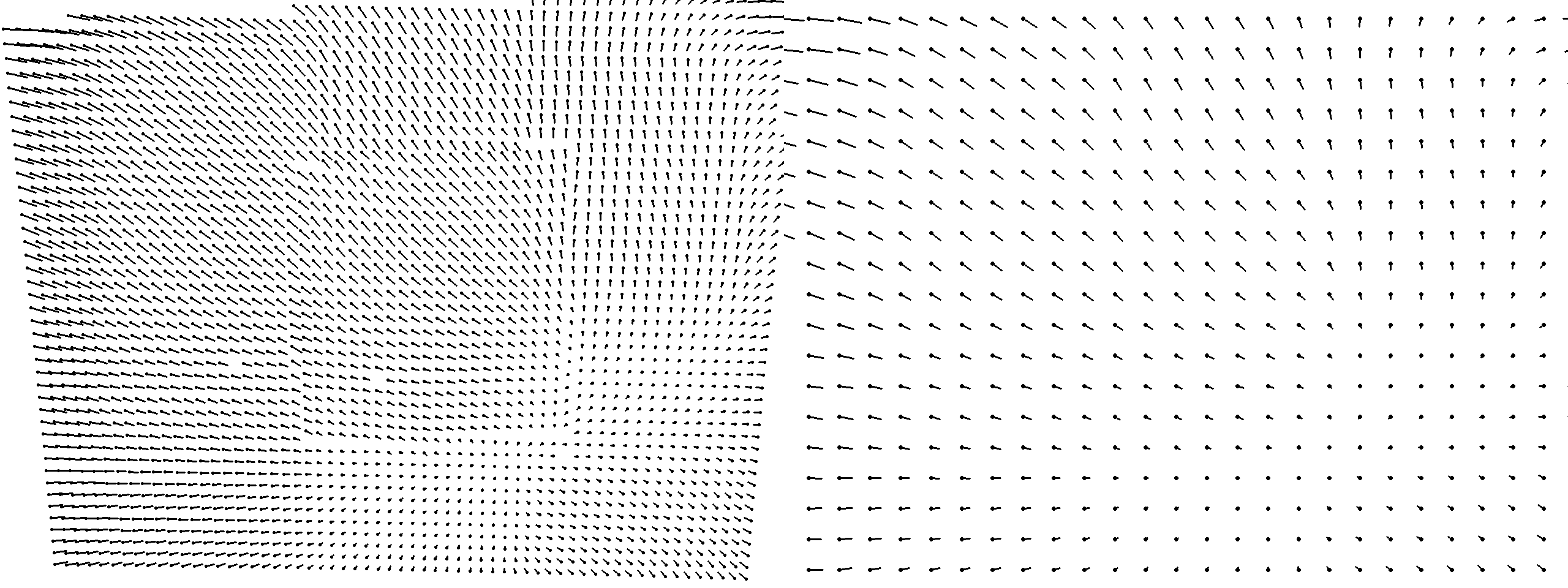}\\
  \caption{Example of selected measured projector distortion vectors (left) and averaged grid representation (right); vectors are always scaled by factor ten.}
  \label{fig:fig3}
\end{figure}

In the final step the vectors $d_{i}$ are averaged to certain predefined grid points in the projector image plane $I_{p}$. The distortion values for all other points of $I_{p}$ being no grid points are interpolated using the grid points. Such, a so called distortion matrix of the projector is obtained (see Fig.~\ref{fig:fig3}). 

\section{Experiments and Results}

\begin{figure}[ht]
  \centering\includegraphics[width=1.0\textwidth]{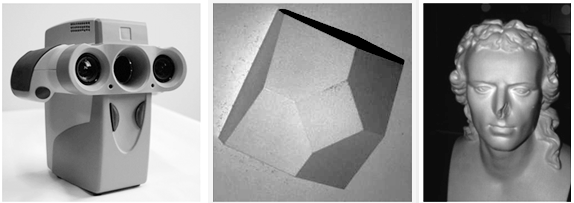}\\
  \caption{Measurement sensor "kolibri Cordless" (left) and measurement objects: stump of a pyramid (middle) and Schiller bust (right).}
  \label{fig:fig4}
\end{figure}

In order to evaluate the described methodology the following experiments were performed. It was used the sensor "kolibri Cordless" (\cite{Munkelt-2007-SPIE}, see Fig.~\ref{fig:fig4}), a mobile, lightweight, hand-held sensor. It is covering a measurement volume of about $280\,mm \times 200\,mm \times 200\,mm$. It has a main triangulation angle between the cameras of about $14^{\circ}$ and between the cameras and the projector of about $7^{\circ}$.

Sequences of phase shifted sinusoidal fringe images with four, eight, and 16 images corresponding to a phase shift of $90^{\circ}$, $45^{\circ}$, and $22.5^{\circ}$ were applied, respectively. The projected width of one fringe period varied between eight and 64 projector pixels, namely 8, 16, 32, and 64 pixels. Additionally, to each sequence a Gray code sequence was added in order to make statements concerning measurement accuracy and completeness in relation to a common reference measurement (GC).

The number of found 3D points using the Gray code sequence was chosen as reference number of "correct" points. However, note that even using Gray code for phase unwrapping, "false" 3D points may be generated, especially at sharp object edges or partly hidden parts of the object surface. Consequently, a completeness level of more than $99 \,\%$ may already include all true 3D points in certain cases.

Two characteristic measurement objects were selected for the experiments concerning completeness of the measurements. The first one was a pyramid stump with a volume of $120\,mm \times 120\,mm \times 40\,mm$ and the second one a Schiller bust with a size of approximately $250\,mm \times 170\,mm \times 190\,mm$ (see Fig.~\ref{fig:fig4}).

\subsection{Measurement Parameters for Evaluation of the New Method}

The following parameters were defined in order to evaluate the new method. Completeness $com$ should be the percentage of the true 3D points in relation to the determined $n$ points using the Gray code measurement. Those 3D points which are not true but recognized as "true" points should be denoted by $fp$ (false positives) as percentage in relation to $n$.

Measurements were performed before (BCI) and after (modes $m_{1}$ and $m_{2}$) calibration improvement by projector distortion correction.
The quality of the calibration should be evaluated by the average length difference $cpd_{mn}$ of the corresponding 3D points obtained by CC-mode and CP-mode application.

\subsection{Results}

Table 1 shows the results for a 16-phase algorithm and a fringe width of 16 pixels for the Schiller bust for both measurement objects before (BCI) and after calibration improvement using modes $m_{1}$ and $m_{2}$ with different values for $thr$. Table 2 shows the results for the pyramid stump. Measurement volume depth was restricted by $200\,mm$ (bust) and $100\,mm$ (stump), respectively.

Before calibration improvement the $cpd_{mn}$ value was about $1.1\,mm$ and after calibration improvement approximately $0.2\,mm$. Epipolar line error $\Delta E$ was smaller than 0.1 pixel.

\begin{table*}[h]
\caption{Percentage results of completeness and false positives (bust)}
\label{tab:Schiller}
\centering
\vspace{+0,3cm}
\begin{tabular}{c | c c c c| c c c c}
\hline
\hline
	thr [rad]	& $0.1\,\pi$ & $0.04\,\pi$	& $0.02\,\pi$	& $0.01\,\pi$	& $0.1\,\pi$ & $0.04\,\pi$ & $0.02\,\pi$ & $0.01\,\pi$ \\
\hline
  Meas.						 & \multicolumn{4}{c|}{$com$ $\left[\%\right]$} & \multicolumn{4}{c}{$fp$ $\left[\%\right]$} \\
\hline
 BCI	& 66.3 & 38.4	& 22.4 & 10.1	& 7.5	& 3.0	& 1.8	& 1.0	\\
 M1 & 98.4 & 95.7 & 91.0 &	73.8 & 2.6 & 2.3 & 2.0 & 1.4 \\
 M2 & 99.1	& 96.3 & 91.5 & 74.2 & 12.7 & 5.5 & 2.7 &	1.5 \\
 GC & \multicolumn{4}{c|} {100}	& \multicolumn{4}{c} {0.2} \\
\hline
\hline
\end{tabular}
\end{table*}

\begin{table*}[h]
\caption{Percentage results of completeness and false positives (stump)}
\label{tab:stump}
\centering
\vspace{+0,3cm}
\begin{tabular}{c| c c c c| c c c c}
\hline
\hline
	thr [rad]	& $0.1\,\pi$ & $0.04\,\pi$	& $0.02\,\pi$	& $0.01\,\pi$	& $0.1\,\pi$ & $0.04\,\pi$ & $0.02\,\pi$ & $0.01\,\pi$ \\
\hline
  Meas.						 & \multicolumn{4}{c|}{$com$ $\left[\%\right]$} & \multicolumn{4}{c}{$fp$ $\left[\%\right]$} \\
\hline
 BCI	& 51.5 & 23.3	& 11.3 & 5.9	& 0.0	& 0.0	& 0.0	& 0.0	\\
 M1 & 99.5 & 98.9 & 98.1 & 86.1 & 0.0 & 0.0 & 0.0 & 0.0 \\
 M2 & 99.5 & 98.9 & 98.1 & 86.1 & 0.0 & 0.0 & 0.0 & 0.0 \\
 GC & \multicolumn{4}{c|} {100}	& \multicolumn{4}{c} {0.0} \\
\hline
\hline
\end{tabular}
\end{table*}

As it can be seen modes $m_{1}$ and $m_{2}$ lead to identical results at the pyramid stump measurement. Regarding the bust measurement mode $m_{1}$ leads to less false positive results but a little bit lower rate of the completeness than mode $m_{2}$. However, it depends on the users priority which mode should be applied. Without calibration improvement (BCI) the results are not sufficient.

The new method has absolutely no influence on the measurement accuracy which only depends on the fringe width and the number of the sinusoidal images in the sequence. At all measurements the standard deviation of the measured 3D surface points was about $11\,\mu m$ for the bust and about $9\,\mu m$ for the stump.

\begin{figure}[ht]
  \centering\includegraphics[width=0.94\textwidth]{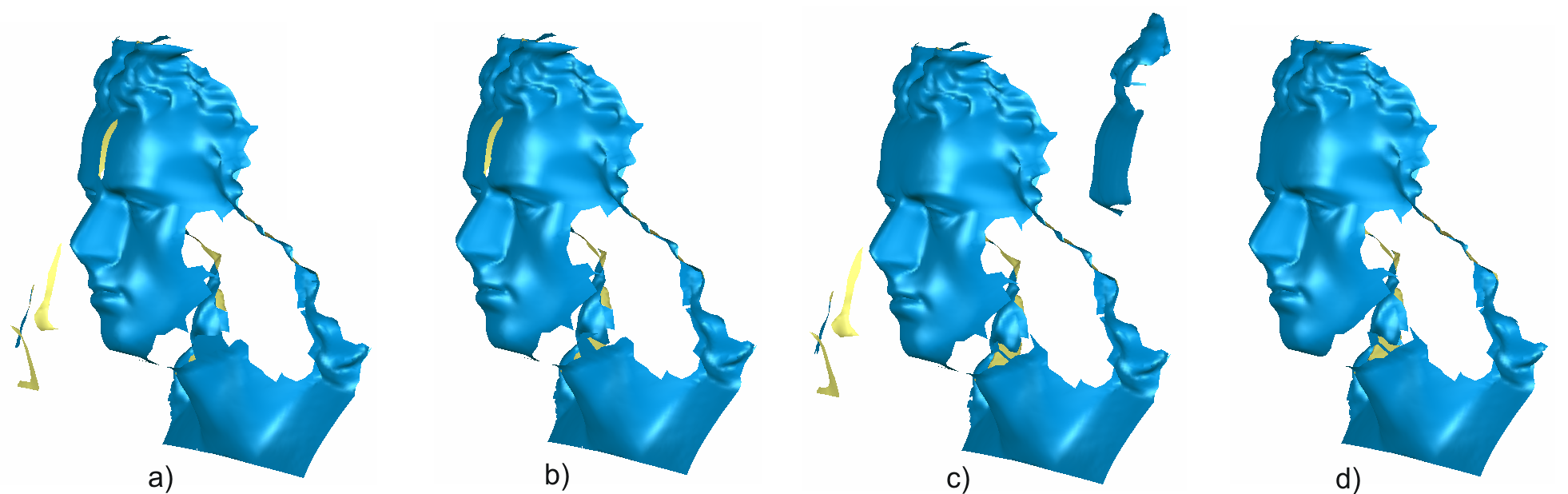}\\
  \caption{Example of a bust measurement: using mode $m_{1}$ without (a), mode $m_{1}$ with (b), mode $m_{2}$ without (c) and mode $m_{2}$ with outlier elimination (d)}
  \label{fig:fig5}
\end{figure}

See Fig.~\ref{fig:fig5} for a typical measurement result of the Schiller bust ($thr = 0.04\,\pi$). Modes $m_{1}$ and $m_{2}$ are compared with and without outlier elimination. Mode $m_{1}$ leads to some missing points on the forehead but has only few false positives. Mode $m_{2}$ leads to a complete result with a relative high number of false positive 3D points. However, because of the big distances to the true 3D points, it is easy to eliminate the false positives manually.

\section{Summary and Outlook}

We presented a new method for geometric phase unwrapping for fringe projection based 3D stereo sensors with reduced projection code. The novelty is the double triangulation for 3D point calculation which leads to almost unique 3D results after considerable calibration improvement of the projector.

Future work should include automatic outlier elimination by 3D point clustering, application of the new method using more scanners, and the improvement of the evaluation of the results. We assume to obtain further improvement of the projector calibration by consideration of 3D distortion effects. If those could be detected and removed, the decision threshold $thr$ could be minimized leading to less calculated false positive 3D points.

\bibliography{refs}

\begin{thebibliography}{10}

\bibitem{CBB-2011-ICIAP}
Christian Br\"{a}uer-Burchardt, Peter K\"{u}hmstedt, and Gunther Notni.
\newblock {Error compensation by sensor re-calibration in fringe projection
  based optical 3D stereo scanners}.
\newblock {\em {Proc. ICIAP}}, pages 363--373, 2011.

\bibitem{CBB-2013-SPIE}
Christian Br\"{a}uer-Burchardt, Peter K\"{u}hmstedt, and Gunther Notni.
\newblock {Phase unwrapping using geometric constraints for high-speed fringe
  projection based 3D measurements}.
\newblock {\em {Proc. SPIE}}, 8789:87890 1--11, 2013.

\bibitem{CBB-2011-Constraints}
Christian Br\"{a}uer-Burchardt, Christoph Munkelt, Matthias Heinze, Peter
  K\"{u}hmstedt, and Gunther Notni.
\newblock {Using Geometric Constraints to Solve the Point Correspondence
  Problem in Fringe Projection Based 3D Measuring Systems}.
\newblock {\em {Proc. ICIAP}}, pages 265--274, 2011.

\bibitem{Herraez-2002-AO}
Miguel~Arevallilo Herr\'{a}ez, David~R. Burton, Michael~J. Lalor, and
  Munther~A. Gdeisat.
\newblock {Fast two-dimensional phaseunwrapping algorithm based on sorting by
  reliability following a noncontinuous path}.
\newblock {\em {Applied Optics}}, 41:7437--7444, 2002.

\bibitem{Ishiyama-2007-AO}
Rui Ishiyama, Shizuo Sakamotom, Johji Tajima, Takayuki Okatani, and Koichiro
  Deguchi.
\newblock {Absolute phase measurements using geometric constraints between
  multiple cameras and projectors}.
\newblock {\em {Applied Optics}}, 46(17):3528--3538, 2007.

\bibitem{Li-2005-OE}
E.B. Li, X.~Peng, J.F. Chicaro, J.Q. Yao, and D.W. Zhang.
\newblock {Multi-frequency and multiple phase-shift sinusoidal fringe
  projection for 3D profilometry}.
\newblock {\em {Optics Express}}, 13:1561--1569, 2005.

\bibitem{Luhmann-2006-Wiley}
Thomas Luhmann, Stuart Robson, Stephen Kyle, and Ian Harley.
\newblock {Close range photogrammetry}.
\newblock {\em {Wiley Whittles Publishing}}, 2006.

\bibitem{Munkelt-2007-SPIE}
Christoph Munkelt, Christian Br\"{a}uer-Burchardt, Peter K\"{u}hmstedt, Ingo
  Schmidt, and Gunther Notni.
\newblock {Cordless hand-held optical 3D sensor}.
\newblock {\em {Proc. SPIE}}, 6618:66180D 1--8, 2007.

\bibitem{Sansoni-1999-AO}
Giovanna Sansoni, Matteo Carocci, and Roberto Rodella.
\newblock {Three-dimensional vision based on a combination of Gray-code and
  phase-shift light projection: Analysis and compensation of the systematic
  errors}.
\newblock {\em {Applied Optics}}, 38(31):6565--6573, 1999.

\bibitem{Schreiber-2000-OE}
Wolfgang Schreiber and Gunther Notni.
\newblock {Theory and arrangements of self-calibrating whole-body
  three-dimensional measurement systems using fringe projection technique}.
\newblock {\em {Optical Engineering}}, 39:159--169, 2000.

\bibitem{VDI-2008-V}
VDI~/ VDE~2634.
\newblock {Optical 3D-measuring systems}.
\newblock {\em {VDI/VDE guidelines, Parts 1-3}}, 2008.

\bibitem{Young-2007-CVPR}
Mark Young, Erik Beeson, James Davis, Szymon Rusinkiewicz, and Ravi
  Ramamoorthi.
\newblock {Viewpoint-coded structured light}.
\newblock {\em {CVPR}}, pages 1--8, 2007.

\bibitem{Zhang-1999-AO}
Hong Zhang, Lalor~Michael J., and David~R. Burton.
\newblock {Spatiotemporal phase unwrapping for the measurement of discontinuous
  objects in dynamic fringe-projection phase-shifting profilometry}.
\newblock {\em {Applied Optics}}, 38(31):3534--3541, 1999.

\bibitem{Zhang-2010-OLE}
Song Zhang.
\newblock {Recent progresses on real-time 3D shape measurement using digital
  fringe projection techniques}.
\newblock {\em {Optics and Lasers in Eng.}}, pages 149--159, 2010.

\end{thebibliography}
\end{document}